\documentclass{amsart}
\usepackage{amstext}

\theoremstyle{definition}

\theoremstyle{remark}

\numberwithin{equation}{section}

\begin{document}
\title{Counterexample to the Hodge conjecture}
\author{K. H. Kim}
\address{K. H. Kim\\
Mathematics Research Group\\
Alabama State University, Montgomery, AL 36101-0271, U.S.A.\\
and Korean Academy of Science and Technology}
\email{khkim@alasu.edu}
\author{F. W. Roush}
\address{F. W. Roush\\
Mathematics Research Group\\
Alabama State University, Montgomery, AL 36101-0271, U.S.A.}
\email{froush@alasu.edu}

\date{}


%\thanks{The authors were partially supported by  NSF Grants
%DMS 9024813 and DMS 9405004. }

\begin{abstract}
  We construct a K3 surface whose transcendental lattice has a self-isomorphism  
which is not a linear combination of self-isomorphisms over $\mathbb{Q}$ which preserve cup products 
up to nonzero multiples.  Products of it with itself give candidates for counterexamples
to the Hodge conjecture which may be of interest.
\end{abstract}
\maketitle
\section{Introduction}

The Hodge conjecture is that every cohomology class over the rational
numbers of a nonsingular
projective algebraic variety $V^n$ and has Hodge filtration $p,p$ is 
represented
by a rational linear combination of classes of subvarieties of complex 
dimension $p$.
Finding such subvarieties is important in many settings such as the
study of fibrations whose fibres are abelian varieties, and of sections
of those fibrations.  There are more general categorical implications, such
as the question of a semisimple category of motives, and the role that 
ordinary
homology plays in stable $A^1$ homotopy \cite{[V1]},\cite{[V4]}, \cite{[SV]} 
and
related theories \cite{[Bl]}.  The definition of the
Hodge structure of the complex cohomology of algebraic
varieties will not be given here, see for instance \cite{[GH]}, 
\cite{[Voi2]}.
This conjecture is a theorem of Lefschetz in complex dimensions $1,n-1$.
Griffiths \cite{[G1]}, \cite{[G2]}, \cite{[G3]} proved that there are 
examples of algebraic cycles which are equivalent
in homology but are not algebraically equivalent; a more recent extension of this 
result is
\cite{[N1]}, and this shows that a kind of counterpart to the Hodge 
conjecture
is false.

In general the Hodge conjecture is unknown even for abelian varieties, that
is, algebraic varieties which are tori differentiably, admitting the group
operation: a discussion
is given in the Appendix by Brent Gordon to \cite{[Lew]}.  In particular it 
is
unknown in general for what are called abelian varieties of Weil type, 
though
Chad Schoen \cite{[Sch]} proved the Hodge conjecture for a few of them.  The 
Hodge structures
of abelian varieties are closely related \cite{[vanG]} to the Hodge 
structures
of a famous class of algebraic surfaces: K3 surfaces are complex structures
differentiably equivalent to a generic quartic hypersurface in $CP^3$.  Yuri Zarhin 
\cite{[Zar]}
proved a number of basic results pertaining to the Hodge structures
of K3 surfaces.  Jordan Rizov \cite{[Riz]} gives a description
of topics related to complex multiplication for the Hodge structures of K3 
surfaces, and proves a fundamental theorem.
A general exposition on complex surfaces including the basics on K3 surfaces 
is
\cite{[VdV]}

Shigeru Mukai \cite{[Mu2]} proved that any Hodge isometry over $\mathbb{Q}$ for K3 
surfaces $X,Y$ is induced by an algebraic cycle $Z\subset X\times Y$.  This was used by Jos\'{e} Mar\'i \cite{[Ma]} to 
prove the Hodge
conjecture for certain products of K3 surfaces.  For our example surfaces related to K3
surfaces are constructed such that the Hodge structure has a linear equivalence to itself over $\mathbb{Q}$,
which is not a linear combination of self-maps that
preserve the cup product quadratic forms.  If there were a suitable partial converse
to Mukai's result in this special situation, then that would give a counterexample
to the Hodge conjecture.  In the general case, a K3 Hodge 
structure will determine the cup products up to rational multiples, but there are 
exceptions.

This present method is not related to our earlier remarks \cite{[KR]}.

\section{ Background on Hodge structures, K3 surfaces, and cohomology 
transfer; polarizations}

Much of this material will follow the references by Gordon \cite{[Lew]} and Rizov \cite{[Riz]}.

Definition.  A Riemann form on a complex torus $V/L$ is a nondegenerate, 
skew-symmetric
$\mathbb{R}$-bilinear form $E:L\times L,V\times V\rightarrow Z, \mathbb{R}$ such that 
$E(x,iy)$ is symmetric
and positive definite.

Complex tori with a Riemann form are abelian varieties.

Definition.  A rational Hodge structure of weight $n$ on a  $\mathbb{Q}$-vector space 
$V$ consists of an
decreasing filtration $F^r(V_{\mathbb{C}})$ of the complexification $V_{\mathbb{C}})$ of $V$ such 
that
$F^r(V_{\mathbb{C}})\oplus \overline F^{n-r+1}(V_{\mathbb{C}}) =V_{\mathbb{C}}  $, $F^0 
(V_{\mathbb{C}})=V_C,F^{n+1}(V_{\mathbb{C}})=0$.

This is equivalent to having a Hodge decomposition $V_{\mathbb{C}}=\oplus 
V_{p,q}=F_p(V_{\mathbb{C}})\cap \overline{F_q(V_{\mathbb{C}} )}$
and is equivalent also to a homomorphism $h_1:U(1)\rightarrow GL(V_{\mathbb{R}})$ of 
real algebraic groups,
given by multiplication by $u^{q-p}$ on $V_{p,q}$, or a homomorphism 
$h:\mathbb{C}^*\rightarrow GL(V_{\mathbb{R}})$ of
real algebraic groups given by multiplication by $z^{-p}\overline{z}^{-q}$.

Definition.  A morphism of rational Hodge structures is a rational linear 
map preserving filtration
(equivalently, any of the other structures in the last paragraph).

Definition.  A polarization of a Hodge structure $(V,h)$ is a morphism $b: 
V\otimes V \rightarrow \mathbb{Q}(-n)$
which as a bilinear form satisfies $b(x,iy), x,y\in V_R$ is symmetric and 
positive definite; $\mathbb{Q}(-n)$
is the weight $2n$ Hodge structure on $\mathbb{Q}$ whose only non-zero summand is $\mathbb{Q}^{n,n}$.

Likewise for an abelian variety, a polarization is an isogeny to the dual 
variety; it is said to
be principal if it is an isomorphism of abelian varieties.  Polarizations of 
abelian varieties are
not unique but the first homology of curves do have canonical polarizations, 
and those polarizations
determine the curves uniquely.  There exist polarizations of the cohomology of
all nonsingular projective algebraic varieties, but again, they may not be unique.
The standard polarizations on primitive cohomology have the form $b(x,y)$ is the evaluation of $xyw$
on the fundamental class of the manifold, for some $w$.

Polarizations exist naturally on tensor products of polarized Hodge structures, and
on substructures of polarized Hodge structures.  

We will be considering situations in which an $H^2$ summand is spanned by cup products from
$H^1$ and the automorphisms of the $H^1$ Hodge structure form a Clifford algebra which
preserves the $H^2$ polarization.

Unless stated otherwise, coefficients in topological homology and cohomology will be
$\mathbb{Q}$.

We will relate these K3 surfaces to the Kuga-Satake construction as treated by van Geemen
\cite{[vanG]},\cite{[KuSa]}.  This construction shows that every polarized weight 2 Hodge structure with
a single holomorphic 2 form which has suitable signature
arises as a Hodge structure from a 1 dimensional polarized Hodge structure.
The 1-dimensional Hodge structure gives a certain abelian variety and
the 2-dimensional Hodge structure will then be a summand on the
2-dimensional cohomology of the abelian variety.

This construction goes as follows.  Start with a vector space $V$ with a bilinear
form $b$ giving the polarization.  Take the Clifford algebra $C$ of this form,
and its even degree part $C^+$.  Let $f_1+if_2$ be a holomorphic 2 form.
Then multiplication in the Clifford algebra by $f_1f_2$ suitably normalized
can be taken as multiplication by $i$ in a complex structure; this extends linearly
to an action of $a+bi$ and hence a Hodge structure of weight 1.  

The Clifford algebra has an anti-involution $i_a$ which reverses the order of
multiplication in monomials from $V$ and a linear map $Tr:C\rightarrow \mathbb{Q}$, the
trace of right multiplication.  Let $e_1,e_2$ be two basis elements for $V$, in a 
basis such that the quadratic form is diagonalized, so that $e_1e_2=-e_2e_1$, and such that the signs of
$e_1^2,e_2^2$ are negative. 
Then a polarization of the Hodge structure on $C^+$ is given by $Tr(\pm e_1e_2 i_a(v) w)$,
for some choice of sign.

This form $E$ allows $C$ to be identified with its dual, and therefore
the natural isomorphism $C^+\otimes C^{+*}\rightarrow End(C^+,C^+)$ gives
an identification of $End(C^+,C^+)$ and $C^+\otimes C^+$.
Embedded in $C^+\otimes
C^+$ are copies of the original Hodge structure on $V$: for a fixed invertible
element $e\in E$ there is a mapping $V \rightarrow End(C^+,C^+)$ which
for $v\in V$ sends $x\in C^+$ to $vxe$.  This map $V\rightarrow C^+\otimes C^+$
commutes with multiplication by the complex numbers and thus gives a sub-Hodge
structure.

Note that for weight 2 Hodge structures, the action of $i$ corresponds to the
square of the action for weight 1 Hodge structures. 

Definition.  The Hodge group of $V$ is the smallest algebraic subgroup of 
$GL(V)$ defined over $\mathbb{Q}$,
containing $h(U(1)).$  The smallest algebraic subgroup containing $h(\mathbb{C}^*)$ 
is the Mumford-Tate
group.

The Mumford-Tate group consists of scalar multiples of elements in the Hodge 
group.  Its representations
determine those aspects of Hodge structures which are important for the 
Hodge conjecture; on any tensor
product of V and its dual, the rational sub-Hodge structures are those 
subspaces which are subrepresentations
of the representation of the Mumford-Tate group.  For a polarizable Hodge 
structure, both
groups will be reductive.   For Hodge structures occurring as the homology 
of a nonsingular algebraic
variety in some dimension, the Hodge group will
preserve a bilinear form obtained from cup products up to scalar multiples, 
which makes it
something like an orthogonal or symplectic group times scalars.

Algebraic representations of a reductive group
are semisimple (even though reductive algebraic groups are not called 
semisimple
as algebraic groups) for instance, \cite{[Milne]}.  

K3 surfaces are 4 dimensional topologically and simply connected, their 
middle cohomology groups are a sum
of 22 copies of $\mathbb{Z}$, the ranks of $H^{0,2},H^{1,1},H^{2,0}$ are 1,20,1.

%Let $A_f$ denote the finite adeles, rational multiples of
%elements of the profinite completion $\hat{Z}$ of $Z$.

A quadratic lattice is a lattice provided with a bilinear form.  $U$ is the 
hyperbolic plane, i.e. it has basis
$e,f$ with $\langle e, e\rangle=\langle f,f\rangle =0,\langle e,f \rangle 
=1$.
$E_8$ is the positive quadratic lattice corresponding to the Dynkin diagram
of type $E_8$.  Then the 2-dimensional cohomology of a K3 surface has type a 
sum of 3 copies of $U$
and 2 of $-E_8$.  This lattice and its form are denoted $L_0,\psi$.  
$L_{2d},\psi_{2d}$ are the sublattices
spanned by $e_1+df_1$ in the first $U$ summand, together with all other 
summands.

$A_X$ is the subgroup of $L_0$, for a K3 surface, generated by Chern 
classes of line bundles,
equivalently by Hodge (1,1)-cycles or algebraic cycles,
and its orthogonal complement is called $T_X$, the transcendental lattice.

Definition. A polarized Hodge structure of degree $d$ of K3 type is a 
homomorphism $h: \mathbb{C}^* \rightarrow
SO(V,\psi_{\mathbb{R}} )$ such that $V,\psi$ as an orthogonal space is equivalent to 
$(V_{2d},\psi_{2d})$,
$\psi$ is a polarization of this Hodge structure, $h$ defines a rational 
Hodge structure of type
$(-1,1),(0,0),(1,-1)$, and the Hodge numbers (ranks) are $1,19,1$.

In the case of $H^2$ polarizations of surfaces and transcendental lattices 
of K3 surfaces,
polarizations of Hodge structures are given by the quadratic forms given by 
cup products; the
different signs for the cup product form vs. the polarization correspond to the different $C^*$ 
action on
$H^{2,0},H^{0,2}$.

Zarhin proved that for Hodge structures of K3 surfaces, the complex 
multiplication, that is, the
algebra $E$ of endomorphisms of the Hodge structure on the transcendental lattice, must be either a 
totally real field or a
quadratic imaginary extension of such a field.  He also proved that the 
Hodge group is
the special orthogonal or unitary group of a form arising from cup products. 
  These
also imply irreducibility of the Hodge structure of the transcendental 
lattice. Thus when any Hodge structure of a nonsingular
projective
algebraic variety is decomposed as a sum of indecomposable Hodge structures, each component 
will
either be isomorphic to that transcendental lattice or have no nonzero homomorphisms
to or from it.

The global Torelli theorem for K3 surfaces was first proved by I. 
Piatetskii-Shapiro and
I. Shafarevich \cite{[PS]}.
Friedman \cite{[Fri]} states the global Torelli theorem, which completely 
determines the set of
possible polarized Hodge structures including cup products over $Z$, for all 
K3 surfaces,
in slightly different notation.  In effect pairs consisting
of a K3 surface and a polarization, a choice of primitive, numerically effective line bundle, are 
in one-to-one correspondence
with lines 
$$\mathbb{C}v\subset L_{2k} \otimes \mathbb{C} \owns v^2=0,v\overline{v}>0$$
modulo automorphisms of $L_0$ over $Z$.   The lines can be identified the 
space of holomorphic
2 forms.  Complex multiplication by a subring of a field $E$ will correspond
to an embedding of $E$ in $\mathbb{C}$ called the reflex field,
acting multiplicatively on the vector $v$,
such there is a compatible ring of additive endomorphisms of $\Lambda$.
A version of global Torelli due to D. Burns and M. Rapaport \cite{[BR]}
also describes mappings among K3 surfaces.

We also require the fact that any birational map of surfaces can be realized
as a regular map by blowups of points each of which introduces only new 
algebraic cycles
in the middle dimensional cohomology \cite{[Hart]} ,\cite{[Hiro]}.  Such blowups
map to the original space and will not affect transcendental summands of the
Hodge structure, nor products from them into $H^4$ except up to integer nonzero
multiples, nor $H^1$ nor products from $H^1$ into $H^2$ modulo algebraic summands
of the Hodge structure.

\section{Quadratic forms}

Some background on quadratic forms over the algebraic number
fields is also needed: By the Hasse-Minkowski
theorem,  two quadratic
forms over an algebraic number field $K$ are isomorphic
if and only if they are isomorphic over every completion of $K$.
Every form can be diagonalized over the $K$, that is, it can be
written in the form $$\sum_{i=1}^n c_ix_i^2.$$
A complete set of invariants for isomorphism of quadratic forms over $K$
is given by the dimension, the signature, the discriminant (determinant of
the matrix of the form up to squares) and the Hasse invariant.
The Hasse invariant of a diagonalized form $\sum c_i x_i^2$
is the Brauer group class of a tensor product of quaternion algebras 
$\otimes_{i\le j} (c_i,c_j)$
over the rational numbers, where $(c_i,c_j)$ is a quaternion algebra
whose norm form in terms of a standard basis is
$x_1^2+c_ix_2^2+c_jx_3^2+c_ic_jx_4^2$. It is simplest to compute the
Hasse invariant locally at each prime, where it can be considered as a 
Hilbert symbol, is $\pm 1$, and
then multiply the Hilbert symbols for all pairs
$(c_i,c_j),i\le j$.  There is a reciprocity formula for Hilbert symbols, 
hence Hasse
invariants, so that knowing the Hasse invariant at all
but one prime for a quadratic form determines its value at the remaining 
prime.
The Hasse invariant of a sum of two quadratic forms is the product of their 
Hasse invariants
times a quantity which depends only on their discriminants.

It will suffice to consider and embed a transcendental lattice over $\mathbb{Q}$,
to study the $\mathbb{Q}$ Hodge structure.  That is we take the quadratic form
over $\mathbb{Z}$, over $\mathbb{Q}$ find a summand of this form, take a mapping, and
carry the period vector with the mapping to construct a K3 surface
with a suitable transcendental lattice.
 Within the transcendental lattice, the 
coordinates
of $z$ must be linearly independent over $\mathbb{Q}$, since any linear combination 
which
is $0$ corresponds to a class in $H^{1,1}$ which is a Hodge cycle, and the
transcendental lattice will not contain any such classes.  For background in
algebraic number theory, see \cite{[W]}.

Proposition 3.1.  For any quadratic forms $q_1,q_2$ over $Q$ such that $q_1$ 
embeds in $q_2$ over $\mathbb{R}$
with codimension at least $4$, $q_1$ will embed in $q_2$ over $\mathbb{Q}$.

Proof.  It will suffice to construct a form of degree 4 whose sum gives 
arbitrary invariants at all
but one prime.  This can be done prime by prime, including the real prime.  
By consideration of the invariants of a sum of two quadratic forms,
it is then enough to construct a 4 dimensional form with arbitrary 
discriminant and Hasse invariants,
except at one non-archimedean prime.
By Theorem 73.1 of \cite{[O]} it is possible to do this globally if it is 
possible to do so at each prime,
given the finiteness of the number of primes involved in the Hasse invariant 
and discriminant.
There are 4 possibilities for the local discriminant (up to multiplication
by local units) and Hasse invariant and they can be realized
by forms as follows
$$\langle 1,-1,1,-1 \rangle , \langle 1,-a,\pi,-a\pi \rangle ,\langle 
1,-1,1,\pi \rangle ,
\langle 1,-1,a,\pi \rangle$$
where $\pi$ is the (odd) prime and $a$ is a nonquadratic residue modulo it.
$\Box$

Proposition 3.2.  Let $E$ be a totally real algebraic number field with
basis $1,b_2,\ldots ,b_{n_0}$.  Consider a Hodge structure with
a single holomorphic 2 form having periods
$$v=b_1, b_2, \ldots b_{n_0}, b_1x_2, x_2b_2, \ldots x_2b_{n_0},
b_1x_3,x_3b_2, \ldots , x_{d_0}b_{n_0}.$$
where $x_i$ are linearly independent over $E$.  Consider a  bilinear
form given by a block diagonal matrix $D$ whose blocks are matrices
of elements of $E$ under the regular representation $F_{11},\ldots 
,F_{d_0d_0}$.  For the regular
representation the basis $b_i$ is used, and it is also assumed that
the regular representation consists of symmetric matrices. (1) Consider the effect of 
complex multiplication by a block diagonal matrix whose blocks 
are some $3\times 3$ matrix $A$ of an element of $E$.  This multiplication 
preserves the period vector and its orthogonal complement up to complex
multiples, and is 
therefore an isomorphism of Hodge structures. 
(2) A quadratic form associated to a group element has signature given
by the number of positive and negative Galois conjugates of that
element (its image under all real embeddings of $E$).  (3) The only 
complex multiplications which preserve the cup product form up 
to rational multiples are by those elements whose squares 
are rational numbers.

Proof. The properties of the regular representation directly implies
that multiplication by the matrix $A$ on the period vector has the
effect of multiplying by the corresponding field element $a$ embedded
in the complex numbers:
$$a(b_1,b_2,b_3)=(ab_1,ab_2,ab_3)$$
$$=(\sum b_ja_{j1} ,\sum b_j a_{j2} ,\sum b_j a_{j3})=(b_1,b_2,b_3)A$$

For the orthogonal complement of $z$, we will be looking at
$wDz$; the effect of multiplying by $A$ is to change this
to $wDAz=wD(az)$ and the orthogonality condition for $w$ is not changed.

The quadratic form given by a symmetric matrix has signature given by
its eigenvalues, and these will be the Galois conjugates of the
field element in a regular representation.

If a complex multiplication by $A$ corresponding to field element $a\in E$
sends all cup products to $k$ times themselves,
$k\in \mathbb{Q}$, then that is also true for complex cohomology, and in particular
for the inner product $z\overline{z}$.  But that inner product under
the complex multiplication goes to $az\overline{az}$, so that $a^2=k$.
$\Box$

Proposition 3.3.  There exists a totally real cubic field and elements
elements $f_1,f_2,f_3$ with the following properties.
The field has a regular representation by symmetric matrices.
The algebraic numbers $f_1,f_2$ have exactly one positive Galois 
conjugate each
and $f_3$ have no positive Galois conjugate.
The positive Galois conjugates of $f_1,f_2$ are much
larger than the other numbers and approximately equal, to any
desired degree, e.g. the positive conjugates
exceed 1 and for any $\epsilon >0$, they differ by at most $\epsilon$,
and the other conjugates are at most $\epsilon$ in
absolute value.  The only elements of the field whose squares
are rational numbers are those which are already rational.

Proof. The statement about Galois conjugates
follows by looking at the Minkowski embedding of the
cubic field in $\mathbb{R}^3$, and its having dense image.  Strong approximation
makes this condition compatible with the other order conditions.
Regarding the symmetric representation, it suffices to take a field
generated by an eigenvalue of some symmetric $3\times 3$
matrix whose characteristic polynomial is irreducible
over the rational numbers. The last statement is true in general since
irrational square roots generate quadratic fields which do not embed
in cubic fields.  $\Box$

\section{Construction of K3 surfaces}

Proposition 4.1. There exists a K3 surface $X_b$ having transcendental lattice of 
rank 9 over K, such that there exists an isomorphism of Hodge structures from 
the lattice
to itself, which is not a linear combination of isomorphisms each of which 
sends
the cup product quadratic form over $\mathbb{Q}$ to a rational multiple of itself.

Proof.  The Hodge structure will have complex multiplication by a
totally real cubic field $E$ of the last proposition.
The transcendental lattice $T_X$ over $\mathbb{Q}$
will be a sum of 3 copies of $E$, as regards complex multiplication.
As in the previous propositions we can choose a family of vectors
$v$, where $x_i$ remain unspecified, which admit this complex multiplication
and also admit cup product forms in which $v$ has the same orthogonal
complement, and
which have the correct signature to embed in $L_{2d}$ for any $d$ as
quadratic forms over $\mathbb{Q}$.  Given any such embeddings the natural mappings
between coordinates of $v$ will give a map which sends $v$ to $v$ and its
orthogonal complement under that bilinear form to its orthogonal complement.
This is then a homomorphism of Hodge structures of two K3 surfaces.

The vectors $v$ are obtained by taking solving for $x_2$ to be close to 1,
pure imaginary and transcendental giving a solution for $x_3$ which is to
be a very small real number.  Then the conditions of  
Proposition 3.3
will ensure positivity of $v\overline{v}$.  There will be no larger field
of complex multiplication.

By Proposition 3.2 complex multiplications preserving the cup product form up
to rational multiples must be square roots of rational numbers, and by 
Proposition 3.3 they are already rational, so complex multiplications by
irrational elements of the cubic field will not be linear combinations of 
isomorphisms of the Hodge structure, which must be some complex multiplications,
which preserve cup products. $\Box$

\section{Kuga-Satake correspondence}
At this point we return to the Kuga-Satake correspondence
\cite{[KuSa]},\cite{[KuSa]} which associates
polarized 1-dimensional Hodge structures $G_1$ with polarized
2-dimensional Hodge structures $G_2$ in which the holomorphic 2-forms are 
all multiples
of a single form $f$, and satisfying a signature condition.

It will be shown that if the Kuga-Satake 
construction is applied to the transcendental lattices of the K3 surfaces 
of Proposition 4.1 then the endomorphisms of the 1-dimensional 
Hodge structure are given by right multiplications in the 
Clifford subalgrebras $C^2$ have unique polarization up to
right multiplication, which will not affect the $H^2$ summands.

Our reason for doing this is the hope that algebraic cycles in 
the Kuga-Satake abelian variety $A$ can be studied in terms of algebraic
cycles giving an equivalence between their 2-dimensional summands
and the 2-dimensional summands of K3 surfaces.

To do 
that in
turn it will suffice to show that the Mumford-Tate groups acting on these 
summands
admit no nontrivial (linear) endomorphisms other than those given by right 
multiplication
(which will not affect $H^2$).  To do that it will suffice to prove the same 
for the
algebra of endomorphisms spanned by the Mumford-Tate group.  This will be 
the algebra
of endomorphisms spanned by right-multiplication by $f_1f_2$ taken in terms 
of a
rational basis.  It will suffice to show these generate the even part of the
Clifford algebra considered as left multiplications.

Proposition 5.1.  For specific transcendental lattices of the K3 surfaces of the last 
section,
there are no nontrivial endomorphisms except for right multiplication within 
the
Clifford algebras.

Proof.  This is a computer calculation using van Geemen's 
constructions and a specific example.  The matrix 
$$A=\left( \begin{matrix} 
               1 &1 &1\\
           1 &1 &0\\
           1 &0 &0
\end{matrix} \right)  $$
has irreducible characteristic polynomial $x^3-2x^2-x+1$.
The elements $b_i$ are chosen so that the regular representation
of the field generated by $A$ with basis $b_i$ will give the
symmetric matrix $A$, and hence the entire regular representation
is by symmetric matrices.  Choose $b_3=1$;
the 3rd row of the matrix of $A$ means we should choose $b_1$
as $A$.  Then the first row means $AA$ should be the sum of the
3 basis elements, $A+I$ plus the 2nd basis element, so the latter
should be $A^2-A-I$.

The period vector $f$ will be chosen in as generic a way as possible
and  $f_1,f_2$ will be its real and complex parts.  We will assume some
 real multiple is chosen to give van Geemen's normalization in terms
of $f\overline{f}$; our results hold for multiples of $f$ by any
real number.

So $f$ is $$b_1,b_2,1,x_1b_1,x_1b_2,x_1,x_2b_1,x_2b_2,x_2$$ or in terms
of a basis $e_i$ for the vector space, later to be generators of the Clifford algebra,
$$b_1e_1+b_2e_21+e_3+x_1b_1e_4+x_1b_2e_5+x_1e_6+x_2b_1e_7+x_2b_2e_8+x_2e_9$$  

We take $x_i$ to be approximations to some solution of $f^2=0,f\overline{f}>0$
which are as general as possible.  One choice of elements 
$\phi_i$ of the cubic field to give the diagonal matrices $F_{ii}$
of the quadratic form is $10I-A-A^2, 2I-7A+2A^2, -A-I$.   When the $x_i$ are chosen in a generic way, then their 
real and imaginary parts and products between them will be linearly independent
over the cubic field.  

The algebraic numbers $\phi_i$ which determine the quadratic form 
satisfy only order conditions and will range over some open neighborhood
within $n$-tuples of rational numbers.  By expanding in terms of a basis, one can see that
if there is a choice of $\phi_i$ within the field such that 
the algebras generated by $x_i$ have rank 256, then the rank will be 256 on the complement of some lower dimensional
algebraic variety for each form, and this complement has elements within any nonempty open
set.  So it is enough to find any field elements such that
the algebra has rank 256.  The particular choice of $\phi_i$
in our computer programs was $A,A^2-A-I,I$.

To compute the van Geemen Clifford algebra, now take the product $f_1f_2$
$$(b_1e_1+b_2e_2+e_3+x_{1r}b_1e_4+x_{1r}b_2e_5+x_{1r}e_6+x_{2r}b_1e_7+x_{2r}b_2e_8+x_{2r}e_9)$$
$$(x_{1i}b_1e_4+x_{1i}b_2e_5+x_{1i}e_6+x_{2i}b_1e_7+x_{2i}b_2e_8+x_{2i}e_9)$$ 
A basis for it over the rational numbers will give elements of the algebra
generated by the Hodge group, since it represents multiplication by $i$ in the
Hodge structure.  For such a basis we can expand out in terms of real and
imaginary parts
$x_{1r},x_{2r},x_{1i}x_{2i}$
and then expand the products of the $b_i$ in terms of the basis $I,A,A^2$
for the cubic field.  This gives nine generators for the algebra generated
by the Hodge group.  They are symmetric to the first case in terms of the 
$x$ expansion, that is, the coefficients of $x_{1r}$ in $f_1f_2$, which are

$$(b_1e_1+b_2e_1+e_3)(b_1e_4+b_2e_5+b_3e_6)$$
$$=(Ae_1+(A^2-A-I)e_2+e_3)(Ae_4+(A^2-A-I)e_5+e_6)$$
$$=A^2e_1e_4+A(A^2-A-I)(e_1e_5+e_2e_4)+(A^2-A-I)^2e_2e_5+A(e_1e_6+e_3e_4)$$
$$+(A^2-A-I)(e_2e_6+e_3e_5)+e_3e_6$$
$$=A^2e_1e_4+(A^2-I)(e_1e_5+e_2e_4)+(A+I)e_2e_5+A(e_1e_6+e_3e_4)$$
$$+(A^2-A-I)(e_2e_6+e_3e_5)+e_3e_6$$

This yields 3 generators
$$e_1e_4+e_2e_6+e_3e_5+e_1e_5+e_2e_4$$ 
$$e_2e_5+e_1e_6+e_3e_4-e_2e_6-e_3e_5$$
$$-e_1e_5-e_2e_4+e_2e_5-e_2e_6-e_3e_5+e_3e_6.$$
In the Clifford algebra, their products will be given by the quadratic
form using $\phi_i$. This proposition will be established if it 
can be shown that these elements (left multiplications) generate the Clifford
algebra, since then no maps but right multiplications will commute with
them, and right multiplications will have trivial effects on $H^2$.
 
The remaining generators differ from the 3 given ones by replacing
$e_1,e_2,e_3$ by $e_4,e_5,e_6$ and $e_4,e_5,e_6$ by $e_7,e_8,e_9$.

A first program computes the regular representation of $e_i$ by $512\times
512$ matrices for
a general quadratic form; then making use of the fact that these
matrices are very sparse, it computes products of $e_ie_j,i<j$.
Then it adds together the sums listed together to get the regular
representation of the 9 generators of the subalgebra.

A second program computes the top row of all 4 fold products of the 9 generators
in which the last generator is one of the first four, and stores
this as a $512\times 729$ array.  A final program puts these matrices 
modulo $101$ into row echelon form.  The computed rank of these products is 
247. By Prop.4.7 the
Clifford algebra $C^+$ is either a full matrix algebra over a quaternion algebra.
By Prop.2.5
in the Appendix to \cite{[Lew]}, the Hodge and Mumford-Tate groups are reductive
and therefore the subalgebra is semisimple and will be a sum of matrix 
algebras over division rings over fields.  In comparing the algebras it is
enough to consider their dimensions over the complex numbers
and a simple algebra of $16\times 16$ matrices will have
no proper subalgebras of rank 247 which are sums of full matrix algebras. 
The rank being 256 means the algebra 
generated by these elements includes all left multiplications.
Therefore there are no symmetries but right multiplications.
$\Box$

Proposition 5.2.  In the Kuga-Satake construction, begin with
a polarized space $V$ of weight two dimension, constructing a Clifford
algebra $C$, a polarized weight weight Hodge structure identified
with $C^+$, and then an embedding of $V\subset C^+\otimes C^+$.  The polarization
with $V$ inherits from this embedding is a
positive rational multiple of the original Hodge structure.  Thus it is independent
of all choices made in the definition of the Kuga-Satake polarization, and is invariant
under right multiplication by the Clifford algebra.

Proof.  A dense subset of polarized Hodge structures on $V$ will have no complex
multiplication except by $\mathbb{Q}$, and hence a unique polarization, in the following
sense.  The Hodge structures are represented in a way similar to those
for K3 surfaces: take diagonalized quadratic forms and general period vectors
with $zz=0,z\overline{z}>0$.  We can take neighborhoods bounded away from zero,
in which the Kuga-Satake construction is continuous, and for a general solution
$z$ there will be no complex multiplication.  The unique polarization must be
both the original polarization, and with some normalization to be within bounds,
the new. One normalization on a given summand might choose a fixed
nonzero integer cohomology class and assume the polarizations have a certain nonzero integer value on 
it paired with itself. Limits of these will have the same property. $\Box$

\end{document}